\documentclass[letterpaper, paper,11pt]{AAS}		

\usepackage{bm}
\usepackage{amsmath}
\usepackage[utf8]{inputenc}
\usepackage{multirow}
\usepackage{breqn}
\usepackage{subfigure}
\usepackage{placeins}
\usepackage[colorlinks=true, pdfstartview=FitV, linkcolor=black, citecolor= black, urlcolor= black]{hyperref}
\usepackage{overcite}
\usepackage{footnpag}			      	
\usepackage{lscape}             
\usepackage{textcomp}           

\PaperNumber{22-565}

\begin{document}

\title{Low Thrust Trajectory Design Using A Semi-Analytic Approach}

\author{Madhusudan Vijayakumar\thanks{Graduate Student, Department of Aerospace Engineering, Iowa State University, USA.},  
Ossama Abdelkhalik\thanks{Corresponding Author, Associate Professor, Department of Aerospace Engineering, Iowa State University, Ames, Iowa 50011, AIAA senior member.}
}

\maketitle{} 		

\begin{abstract}
Space missions that use low-thrust propulsion technology are becoming increasingly popular since they utilize propellant more efficiently and thus reduce mission costs. However, optimizing continuous-thrust trajectories is complex, time-consuming, and extremely sensitive to initial guesses. Hence, generating approximate trajectories that can be used as reliable initial guesses in trajectory generators is essential. This paper presents a semi-analytic approach for designing planar and three-dimensional trajectories using Hills equations. The spacecraft is assumed to be acted upon by a constant thrust acceleration magnitude. The proposed equations are employed in a Nonlinear Programming Problem (NLP) solver to obtain the thrust directions. Their applicability is tested for various design scenarios like orbit raising, orbit insertion, and rendezvous. The trajectory solutions are then validated as initial guesses in high-fidelity optimal control tools. The usefulness of this method lies in the preliminary stages of low-thrust mission design, where speed and reliability are key.
\end{abstract}

\section{Introduction}

Design of space missions that utilize continuous low thrust propulsion is becoming increasingly popular as they utilize propellant more effectively and ever since has gained much attention in the literature. One of the fundamental task of such a mission is the design of spacecraft trajectory that delivers the spacecraft from a given state to a desired state within a specified amount of time. However, this process is highly challenging, because there are too many trajectory parameters like launch and arrival date, position, velocity, etc which need analyzing. In addition, to this, long duration thrust profiles associated with low-thrust spacecraft must be obtained. The search space associated with these parameters is very large and further compound the complexity of the trajectory design problem. 

At the preliminary stage, thousands, if not millions of possible trajectories are evaluated before obtaining a feasible trajectory for the mission. The generation of guidance trajectory is treated as an optimization problem and usually takes several days or even months to complete. For a general case of low-thrust trajectory design, no analytical closed form solution exists till date \cite{conway2010spacecraft}. The trajectory is designed typically as a boundary value problem using direct or indirect optimization methods \cite{betts1998survey}. The indirect method resorts to calculus of variation while the direct method utilizes nonlinear programming to solve this problem. Direct methods require initial guesses while indirect methods are extremely sensitive to initial guesses and have a very small radius of convergence. Hence direct methods are usually preferred in trajectory design. One of the noticeable works in the direction of direct method was proposed by Jon A. Sims et al. \cite{sims2006implementation}. In this work, the trajectory design problem was divided into multiple discrete segments and solved using an Non-Linear Programming (NLP). Additionally, the branch of shape based methods introduced by Petropolous et al. \cite{Petropoulos2004Longuski} improve the convergence of direct solvers by providing approximate initial trajectory solutions. Early work in this direction include the exponential sinusoidal \cite{Petropoulos2004Longuski} and inverse polynomial method \cite{wall2009shape} in which radical vector of spacecraft is written as a function of the transfer angle. These methods are useful in generating orbit raising trajectories. The shaping pseudoequinoctial \cite{de2006preliminary} and Finite Fourier Series \cite{taheri2012shape} method extended the shape based methods to generate trajectory solutions to a variety of problems like rendezvous, orbit insertion, etc capable of handling thrust constraints.

One of the primary purpose of the approximate trajectory generation methods is to provide a quick evaluation of the search space to look for regions of the design space associated with lower mission costs and avoid lengthy unnecessary calculations. To support this purpose several authors have studied orbit motion and provided analytic solutions for some special cases of orbit motion. Studying the special cases of radial thrust for escape trajectories from circular orbits, Tsien \cite{tsien1953take}, Boltz \cite{boltz1991orbital}, Prussing \cite{prussing1998constant} and Mengali \cite{mengali2009rapid} developed analytical solutions of orbit motion. Megali et al. \cite{mengali2009escape} later extended this study for elliptical orbit departure. Similarly,  Zee \cite{zee1968low}, Boltz \cite{boltz1992orbital}, and Benney \cite{benney1958escape} studied the case of tangential thrust for continuous low thrust trajectories. Furthermore, Gao \cite{gao2005analytic} presented an averaging technique to obtain analytical solution in case of tangential thrust. Although those methods can, to a large extent, simplify the continuous low thrust problem, they are only suitable for some special cases. 

One of the key analytical contribution namely Clohessy–Wiltshire (CW)  equations \cite{clohessy1960terminal} and the Tschauner-Hempel equations \cite{tschauner1965rendezvous} that were developed in the 1960s changed the game of trajectory design. These equations provide linear models for the study of spacecraft relative dynamics in circular and elliptic orbits. For decades, these sets of equations have been the reference model for the design of relative guidance, navigation, and control systems \cite{fehse2003automated}. While the original formulation of the CW equations captures only the rectilinear relative states of the spacecraft, more recent developments by Alfriend et al. \cite{alfriend2009spacecraft} and De Bruijn et al. \cite{righetti2011handling} formulate the CW equations using curvilinear relative states. However, these formulations assume that the spacecraft is not acted upon by any orbital perturbation. Later, Leonard et al. \cite{leonard1987orbital}, Humi et al. \cite{humi2002rendezvous}, Carter et al. \cite{carter2002clohessy}, incorporated the effects of atmospheric drag into these formulations. Bevilacqua et al. \cite{bevilacqua2008rendezvous} demonstrated the extension of this formulation into the domain of continuous low thrust. 

Recent work by Takao et al. \cite{takao2022delta} simplified the state-space formulation of the two-dimensional CW equations with low thrust propulsion as a perturbed acceleration and demonstrated the design of gravity assist interplanetary trajectories. Further, Toshinori et al. \cite{ikenaga2015interplanetary} extended this work to provide a flexible orbit design method for designing low thrust missions. Taking inspiration from the state space implementation of the continuous low thrust CW equations, this paper provides a three-dimensional linear formulation of the Clohessy–Wiltshire equations. The 3D analytic approximation of the CW equations are then used to develop a robust general semi-analytic trajectory design algorithm that is capable of generating both two-dimensional and three-dimensional low thrust trajectories. The applicability of the algorithm is tested for various design scenarios like orbit raising, orbit insertion and rendezvous. The accuracy of the trajectory solutions is evaluated by drawing comparison to numerical propagation. The key feature of this algorithm is in its ability to handle many revolution low thrust trajectories.

This paper is organized as follows. Section II summarizes the low thrust analytic approximation of the Clohessy–Wiltshire equations. Section III presents a brief description for the trajectory design algorithm. Section IV illustrates some numerical examples of both 2D and 3D trajectory design. Section V presents the conclusion of this work.

\section{Low Thrust Analytic Approximation of Hill's Equations}

\begin{figure}[h!]
  \centering
  \begin{minipage}{0.49\textwidth}
    \includegraphics[width=\textwidth]{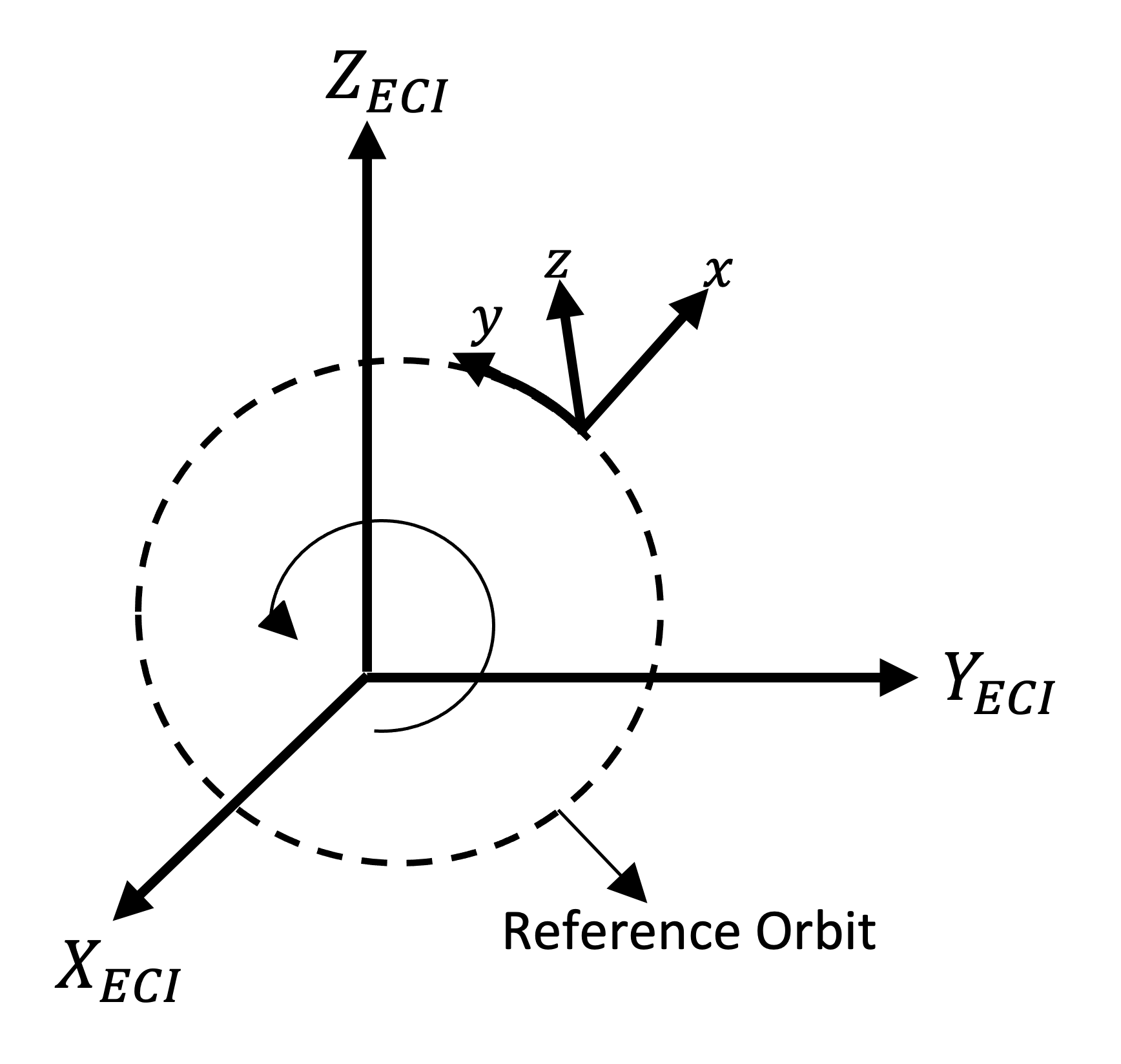}
    \caption{Coordinate Frames}
    \label{fig:CF}
  \end{minipage}
  \hfill
  \begin{minipage}{0.49\textwidth}
    \includegraphics[width=\textwidth]{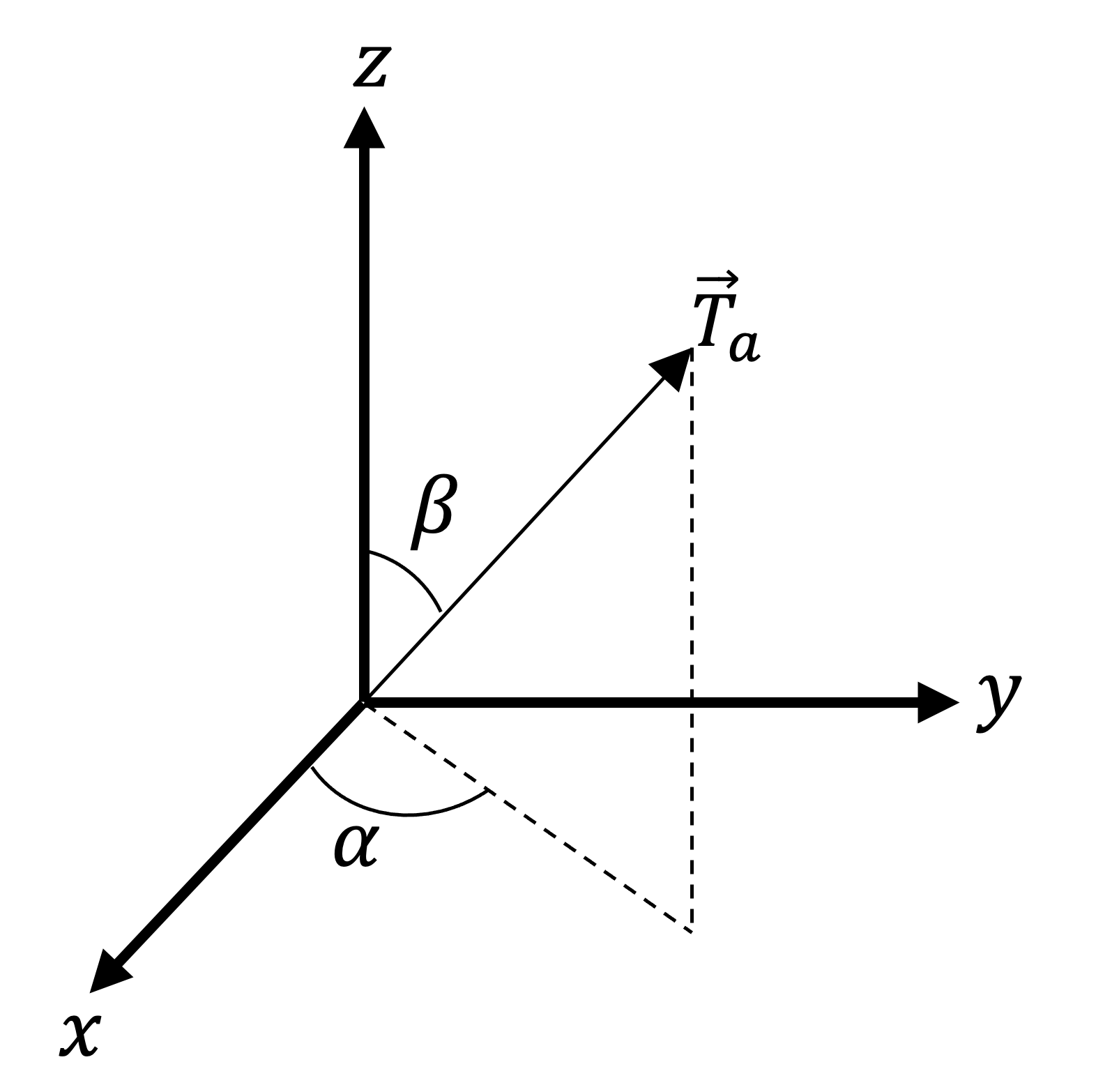}
    \caption{Hills Frame}
    \label{fig:Hills_Frame}
  \end{minipage}
\end{figure}

The Clohessy-Wiltshire equations, also referred as Hill's equations are a set of linear, time invariant differential equations used to describe the motion of a spacecraft relative to a reference orbit. The reference orbit is usually a circle and the corresponding coordinate frame of the relative motion is shown in Fig ~\ref{fig:CF}. The $X-Y-Z$ frame represents the inertial coordinate system and $x-y-z$ represents the Hill's frame. In our case, $X-Y-Z$ represents the Earth Centered Inertial Frame (ECI). The $x$, $y$ and $z$ coordinates of the Hill's frame are the radial, along-track and cross-track displacements of the spacecraft relative to the reference orbit.
The primary goal is to obtain an analytic approximation Consider the two-body Hill's equation as shown in Eq.~\eqref{CW_Eqn}.
\begin{subequations}
\label{CW_Eqn}
\begin{align}
     \ddot{x}-2\,x\,n^2 -2\,\dot{y} \,n &=a_x \label{eq:CW_x} \\
     \ddot{y} +2\,n\,\dot{x} &=a_y \label{eq:CW_y} \\
    \ddot{z} + z\,n^2 &= a_z \label{eq:CW_z}
    \end{align}
\end{subequations}
where, $a_x$, $a_y$ and $a_z$ are the components of the thrust acceleration represented in the Hill's frame. The phase angle of the in-plane thrust acceleration component varies linearly while that of the out-of-plane component is assumed to be constant along each time step as show in Fig ~\ref{fig:Hills_Frame}. Let $\alpha_0$ be the initial phase angle of the in-plane thrust acceleration component and $k$ be the rate at which it rotates in the Hill's frame. Then the thrust acceleration components can be represented as:
\begin{equation}
\begin{aligned}
    a_x &= a\sin(\beta)cos(\alpha_0 - kt) \\
    a_y &= a\sin(\beta)sin(\alpha_0 - kt) \\
    a_z &= a\cos(\beta)
\end{aligned}
\end{equation}
The primary task is to obtain an approximate analytic solution of the  Clohessy-Wiltshire equations for spacecrafts employing continuous low thrust. One way to solve the second-order homogeneous differential equations in Eq.~\eqref{CW_Eqn}, is to use Laplace transform. The out-of-plane motion in Eq.~\eqref{eq:CW_z} is uncoupled from the in-plane motion i.e. only the $z$ component appears and is the simplest among them to solve. This equation resembles the equation of an undamped harmonic oscillator with a forcing (or input) function on the right-hand side. Taking the Laplace transform of \label{CW_z} and rewriting the equation we get:
\begin{equation}
\label{Z_laplace}
    Z\left(s\right)=\frac{{\dot{z} }_0 +s\,z_0 +\frac{a\,\sin{\beta}}{s}}{n^2 +s^2 }
\end{equation}
where, $s$ is a complex number, $z_0$ and ${\dot{z}}_0$ are the z-component of the velocity vector at $t = 0$. Taking the inverse Laplace transform of Eq.~\eqref{Z_laplace} we obtain:
\begin{equation}
    \label{zt_eqn}
    z(t) = \frac{n^2 \,z_0 \,\mathrm{cos}\left(n\,t\right)-a\,\mathrm{cos}\left(n\,t\right)\,\mathrm{sin}\left(\beta \right)+n\,{\dot{z} }_0 \,\mathrm{sin}\left(n\,t\right)}{n^2 }+\frac{a\,\mathrm{sin}\left(\beta \right)}{n^2 }
\end{equation}
Differentiating the above equation we get: 
\begin{equation}
\label{zdott_eqn}
    \dot{z}(t) = {\dot{z} }_0 \,\mathrm{cos}\left(n\,t\right)-n\,z_0 \,\mathrm{sin}\left(n\,t\right)+\frac{a\,\mathrm{sin}\left(n\,t\right)\,\mathrm{sin}\left(\beta \right)}{n}
\end{equation}
From Eq.~\eqref{zt_eqn} and Eq.~\eqref{zdott_eqn} notice that the out-of-plane motion is only influenced by the initial position $z_0$, initial velocity $\dot{z_0}$ and time $t$. Thus, we can write, 
\begin{equation}
\begin{aligned}
\label{eq:y_function}
    z(t) = f_z (z_0,\dot{z}_0,\beta,t) \\
    \dot{z}(t) = f_{\dot{z}} (z_0,\dot{z}_0,\beta,t)
    \end{aligned}
\end{equation}
On the contrary, Eq.~\eqref{eq:CW_x} and ~\eqref{eq:CW_y} highlights that the planar motion of the spacecraft along $x$ and $y$ directions are coupled to each other. The closed form analytic solution for the in-plane position and velocity are obtained as follows. Start from Eq.~\eqref{eq:CW_x} and differentiate it to get:
\begin{equation}
\dddot{x} = 2\,n\,\ddot{y}+3\,n^2 \dot{x} + a\,k\,\mathrm{sin}\left(\alpha_0 -k\,t\right)\,\mathrm{cos}\left(\beta \right)
\end{equation}
Now, substitute $\ddot{y} = -2\,n\,\dot{x} + a\,\mathrm{cos}(\beta)\,\mathrm{sin}\left(\alpha_0 -k \,t\right)$ in the $\dddot{x}$ to get:
\begin{equation}
\dddot{x} = -n^2 \,\dot{x} + a\,\mathrm{sin}\left(\alpha_0 -k\,t\right)\,\mathrm{cos}\left(\beta \right)\,{\left(k+2\,n\right)};
\end{equation}
Take the Laplace transform of the $\ddot{x}$ to get,
\begin{equation}
X\left(s\right)=\frac{{\ddot{x} }_0 +s\,{\dot{x} }_0 +n^2 \,x_0 +s^2 \,x_0 -\frac{a\,\mathrm{cos}\left(\beta \right)\,{\left(k\,\mathrm{cos}\left(\alpha_0 \right)-s\,\mathrm{sin}\left(\alpha_0 \right)\right)}\,{\left(k+2\,n\right)}}{k^2 +s^2 }}{n^2 \,s+s^3 }
\end{equation}
Now, taking the inverse Laplace transform transform and substituting $\ddot{x}_0$ from Eq.~\eqref{eq:CW_x} we get an expression for $x(t)$. Further differentiating the $x(t)$ equation, we can arrive an expression for the velocity $\dot{x}(t)$. Notice that both $x(t)$ and $\dot{x}(t)$ are only a function of the planar components of the initial states as shown in Eq.~\eqref{eq:x_function}. Detailed expressions for the position and velocity $x(t)$ and $\dot{x}(t)$ is presented in the Appendix.
\begin{equation}
\begin{aligned}
\label{eq:x_function}
    x(t) = f_x (x_0,\dot{x}_0,\alpha_0,\beta,k,t) \\
    \dot{x}(t) = f_{\dot{x}} (x_0,\dot{x}_0,\alpha_0,\beta,k,t)
    \end{aligned}
\end{equation}
By substituting the expression for $\dot{x}(t)$ from Eq.~\eqref{eq:x_function} into the Hill's equation in Eq.~\eqref{eq:CW_y}, one can easily obtain expressions for $\dot{y}(t)$ and $y(t)$ by successive integration with respect to time. This would result in both $y(t)$ and $\dot{y}(t)$ expressed as functions of the planar components of the initial states as shown in Eq.~\eqref{eq:y_function}. Detailed expressions for the position and velocity $y(t)$ and $\dot{y}(t)$ is presented in the Appendix.
\begin{equation}
\begin{aligned}
\label{eq:y_function}
    y(t) = f_y (x_0,\dot{x}_0,\dot{y}_0,\alpha_0,\beta,k,t) \\
    \dot{y}(t) = f_{\dot{y}} (x_0,\dot{x}_0,\dot{y}_0,\alpha_0,\beta,k,t)
    \end{aligned}
\end{equation}
Since the closed-form analytic solution of the spacecraft states is obtained by the linearization of Hill's equations, it is important to test the accuracy of the solution. One way to estimate the accuracy preserved by this method is to examine the deviation in the spacecraft states compared to its intended states. This can be easily done by numerically propagating the states of the spacecraft using the two-body equations of motion and comparing the results with the analytic approximation. This experiment was conducted for a series of orbits ranging from LEO to GEO with varying initial conditions. However, only two of the cases are demonstrated here since the general trend remained the same. 

\begin{figure}[h!]
  \centering
  \begin{minipage}{0.49\textwidth}
    \includegraphics[width=\textwidth]{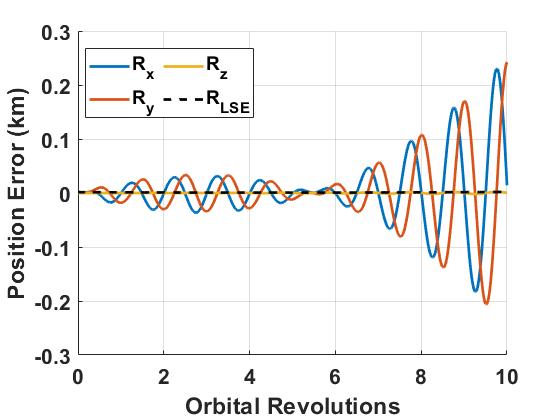}
    \caption{Hill's vs Keplerian Motion - Circular}
    \label{fig:Err_LEO}
  \end{minipage}
  \hfill
  \begin{minipage}{0.49\textwidth}
    \includegraphics[width=\textwidth]{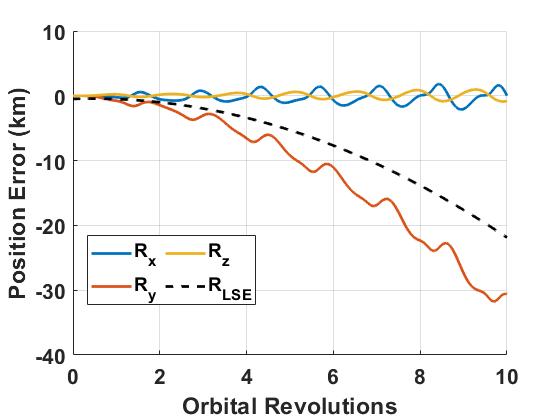}
    \caption{Hill's vs Keplerian Motion - Eccentric}
    \label{fig:Err_HEO}
  \end{minipage}
\end{figure}

In the first case, the spacecraft starts an initial circular low Earth orbit (a = $6678$ \, km, e = $0$) with a period of ~$90$ min. It was assumed that the spacecraft was acted upon by a constant magnitude thrust acceleration of $9e-8$\,km/s$^2$. The initial thrust steering angles were kept at $\alpha_0 = 90^\circ$ and $\beta = 45^\circ$. Figure~\ref{fig:Err_LEO} shows the corresponding positional difference between the analytical approximation and the two-body numerical propagation of the spacecraft. Solid curves are used to indicate the positional difference along each component in the ECI space, and the dashed curve represents the least square fit of the positional error.
For circular orbits, the analytical approximation closely matches the Keplerian motion. From Fig~\ref{fig:Err_LEO}, one can notice that the maximum positional difference is about 0.3 km for the first ten revolutions of the spacecraft. This is the consequence of the circular orbit restriction assumed in Hill's equations. To further examine the accuracy of the analytic approximation, the motion of the spacecraft starting from an eccentric orbit is studied. In this case, the spacecraft starts an initial eccentric low Earth orbit (a = $8164$ \, km, e = $0.17$) with a period of ~$120$ min. It was assumed that the spacecraft was acted upon by a constant magnitude thrust acceleration of $9e-8$\,km/s$^2$. The initial thrust steering angles were kept at $\alpha_0 = 90^\circ$ and $\beta = 45^\circ$. Figure~\ref{fig:Err_LEO} shows the corresponding positional difference between the analytical approximation and the two-body numerical propagation of the spacecraft. Solid curves are used to indicate the positional difference along the each component in the ECI space and the dashed curve represents the least square fit of the positional error. Notice there is a significant decrease in the prediction accuracy of the analytic approximation compared to Keplerian motion. This is not unexpected as the analytic approximation assumes circular motion and the $\omega t$ term experiences large variations in eccentric orbits. One way to overcome this is to restart Hill's propagation at the end of each revolution assuming circular orbits. However, the $\omega t$ term would still show maximum deviation from the mean value for each revolution. This effect is exaggerated in lower altitudes due to shorter orbital periods. Hence, eccentric orbits in LEO are the worst place to use the analytic approximation.

Additional observations from the series of these experiments are summarized as follows. Firstly, error in the decoupled out-of-plane motion are very small compared to the error in the in-plane motion. This can be associated with the fact that the sources for error in the $z$-direction is half as that of the $x-y$ plane. Secondly, irrespective of the initial orbit, the error in the analytical approximation is negligible for propagation times up to one orbital period of the reference orbit. Beyond this, the errors associated with the analytical approximation increases due to reasons mentioned before.

\section{Trajectory Design}
This section provides details on the trajectory design method employed using the analytic approximation of Hill's equations for low thrust trajectories. The error from the circular orbit restriction in the analytic approximation is evident looking at the positional errors studied in the previous section. To overcome this issue, the entire trajectory design is divided into 'm' segments. Each segment of the trajectory is then approximated using the analytic solution. Given the initial conditions (initial states), the trajectory design algorithm iterates over the time of flight (ToF), the thrust steering angles ($\alpha$ and $\beta$) and the rate of change of thrust direction ($k$) along each segment to achieve specific target conditions. The target orbits discussed in this work correspond to either orbit raising, orbit insertion or rendezvous problems. The continuity of the trajectory is enforced by using the final states of the earlier segment as the initial condition of the later segment. 

The entire trajectory design problem is formulated as an optimization problem where the objective is to match the states of the spacecraft in the target orbit while minimizing the total $\Delta V$. The design and objective space of the optimization problem is manifested differently based on the problem being solved. For a rendezvous problem, the design variables are the time of flight (ToF), the thrust steering angles ($\alpha$ and $\beta$) and the rate of change of thrust direction ($k$) along each segment. The target states of the spacecraft are enforced as an equality constraint while the objective is to minimize the total $\Delta V$. The formulation is shown as follows:
\begin{eqnarray}
    \textbf{Minimize  } \,\, J & = & \Delta V   \\ 
    \textbf{Subject to  } \, \, C&:& \begin{bmatrix} a_f \\ e_f \\ i_f \\ \omega_f \\ \Omega_f \\ \theta_f \end{bmatrix} = \begin{bmatrix} a_{t} \\ e_{t} \\ i_{t} \\ \omega_{t} \\ \Omega_{t} \\ \theta_{t} \end{bmatrix} \\
    \textbf{Design Variables  } \, \, X & = & [ToF, \, \, \alpha_0, \, \, \beta_i, \, \, k_i]_{1 \times 2m+2} \, \, ,  \, \, i = 1,2, \cdots m
\end{eqnarray}
where, [$a, e, i, \omega, \Omega, \theta$] are the classical orbital elements. The subscript 'f' represents the final states of the spacecraft at the end of the analytic approximation and the subscript 't' represents the target orbit. The orbit insertion problem is treated as a special case of the rendezvous problem where the true anomaly ($\theta$) is set free while the rest of the formulation remains the same.

For an orbit raising problem, the time of flight becomes a fixed parameter while the objective is to increase the semi-major axis of the final orbit. Additionally, there is no need to enforce any equality conditions. The formulation of the optimization problem is as follows:
\begin{eqnarray}
    \textbf{Maximize  } \,\, J & = & a_f   \\ 
    \textbf{Design Variables  } \, \, X & = & [\alpha_0, \, \, \beta_i, \, \, k_i]_{1 \times 2m+1} \, \, ,  \, \, i = 1,2, \cdots m
\end{eqnarray}

\section{Case Studies}
All the calculations for the case studies are carried out in canonical units such that one distance unit (DU) is equal to the radius of the Earth of 6378.14 km, one time unit (TU) is 806.8 seconds. The canonical units were scaled up by a scaling factor (p) for the case studies that start at high Earth Orbits (HEO) and Geo Stationary Orbits (GEO) as specified in their corresponding sections. All the calculations were carried out in Matlab. The in-built function `fmincon' was used as the optimization routine. All cases were simulated using Intel(R) Xeon(R) Quadcore Processor @ 3.50GHz.

\subsection{Case 1: Coplanar Transfer LEO-to-GEO (Orbit Insertion)}\label{Case1} 

A low-thrust trajectory design problem from an initial circular low Earth orbit (LEO) at an altitude $h_i = 300 \,$ km to a final circular Geostationary Orbit (GEO) of altitude $h_f = 35,785\,km$  is studied in this section. This is a typical orbit insertion problem where the spacecraft is acted upon by a fixed magnitude of thrust acceleration while the direction is variable. For this problem, the parameters of the spacecraft's propulsion system are defined as follows: specific impulse Isp = $3300\,s$, initial spacecraft mass = $95 \, kg$, thrust-to-weight ratio = $10^{-4}$ and maximum thrust limit = $0.0920\,N$ and constant thrust acceleration of Ta = $9.81e-7\,m/s^2$.

\begin{figure}[h!]
  \centering
  \begin{minipage}{0.49\textwidth}
    \includegraphics[width=\textwidth]{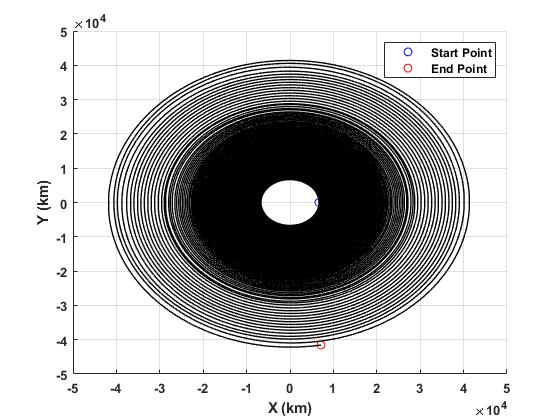}
    \caption{Case 1: Trajectory (ECI Frame)}
    \label{fig:Case1_LEO_GEO_Traj}
  \end{minipage}
  \hfill
  \begin{minipage}{0.49\textwidth}
    \includegraphics[width=\textwidth]{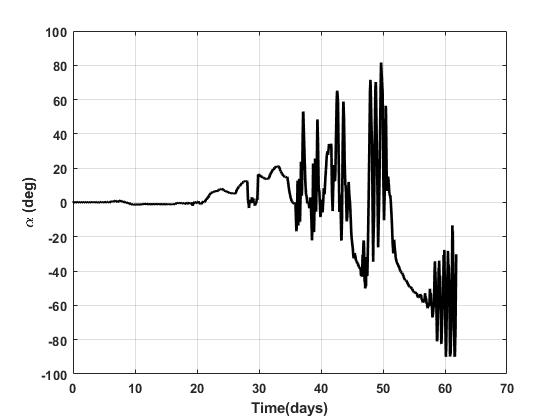}
    \caption{Case 1: History of $\alpha$}
    \label{fig:Case1_LEO_GEO_alpha}
  \end{minipage}
\end{figure}

The trajectory solution provided by the proposed semi-analytic method is shown in Fig~\ref{fig:Case1_LEO_GEO_Traj}. 
The spacecraft takes $61.80\,days$ to transfer from LEO to GEO performing $N_{rev} = 369 \,$ revolutions around the Earth. In doing so, it consumes $14.01\,kg$ of propellant and achieves a $\Delta V = 5.23 \,km/s$. The corresponding thrust profile is an oscillating sinusoid in the $x$ and $y$ direction. The history of the thrust steering angle is as shown in  Fig~\ref{fig:Case1_LEO_GEO_alpha}.
\begin{figure}[h!]
  \centering
  \begin{minipage}{0.49\textwidth}
    \includegraphics[width=\textwidth]{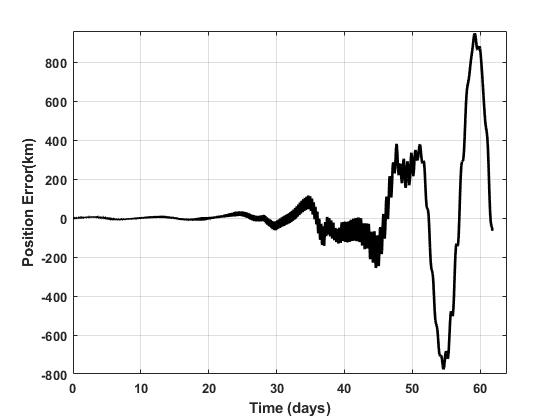}
    \caption{Case 1: Positional Error}
    \label{fig:Case1_LEO_GEO_r_err}
  \end{minipage}
  \hfill
  \begin{minipage}{0.49\textwidth}
    \includegraphics[width=\textwidth]{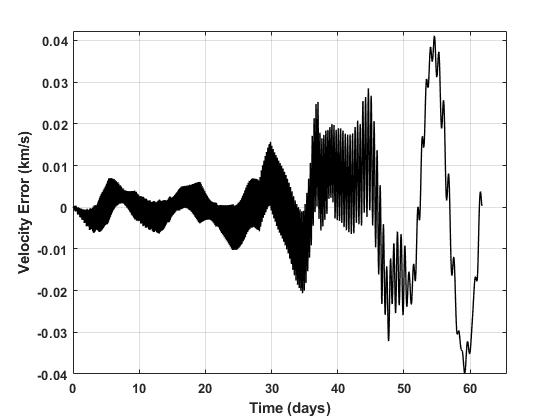}
    \caption{Case 1: Velocity Error}
    \label{fig:Case1_LEO_GEO_v_err}
  \end{minipage}
\end{figure}
Figures~\ref{fig:Case1_LEO_GEO_r_err} and ~\ref{fig:Case1_LEO_GEO_v_err} illustrate the errors associated with the position and velocity of the analytically approximate trajectory compared to numerical propagation. Additionally, Table~\ref{tab:Case1_Performance_Characteristics} provides a detailed list of the obtained solution characteristics.
\begin{table}[h!]
	\fontsize{10}{10}\selectfont
    \caption{Case 1: Performance Characteristics of the Semi-Analytic Solution}
   \label{tab:Case1_Performance_Characteristics}
        \centering 
   \begin{tabular}{c | c } 
      \hline 
        Parameters & Value \\ 
            \hline
        No. of revolutions & 338 \\
        Time of Flight (days) & 45.44 \\ 
        $\Delta V$ (km/s) & 3.85 \\
        Propellant Consumption (kg) & 10.52 \\
      \hline
   \end{tabular}
\end{table}
The final states of the spacecraft achieved by the semi-analytic method and numerical propagation are compared in Table~\ref{tab:Case1_CW_NP_Comp} to demonstrate the level of error between them. From the table it one can notice that the error in matching the target position is about $65$ km and that of target velocity is about $0.005$ km/s. This small error in matching the final conditions obtained by numerical propagation demonstrates the ability of the semi-analytic solution to provide high accuracy solutions.

\begin{table}[hbt!]
\caption{\label{tab:Case1_CW_NP_Comp} Case 1: Semi-Analytic Solution vs Numerical Propagation}
\centering
\begin{tabular}{ccc}
\hline
\hline
        Target orbit Parameters & Semi-Analytic Solution & Numerical Propagation \\ 
        \hline
        Final Position (km) & 42165  & 42100  \\
        Final Velocity (km/s) & 3.0746 & 3.0751\\
        Semi-major axis, $a$ (km) & 42165.18 & 42047.93 \\ 
        Eccentricity, $e$ (km) & 0 & 0.001 \\ 
\hline
\hline
\end{tabular}
\end{table}

\subsection{Case 2: Coplanar Transfer LEO-to-MEO (Rendezvous)}\label{Case2} 

\begin{figure}[h!]
  \centering
  \begin{minipage}{0.49\textwidth}
    \includegraphics[width=\textwidth]{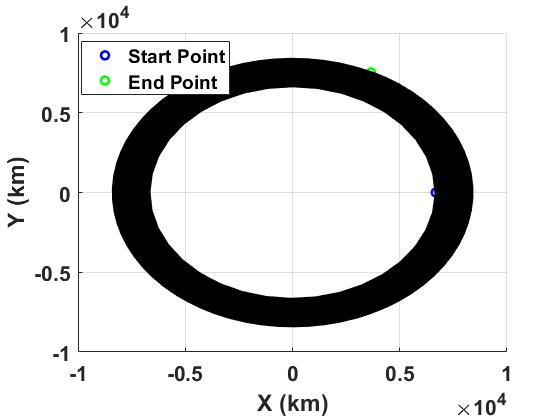}
    \caption{Case 2: Trajectory (ECI Frame)}
    \label{fig:Case2_LEO_MEO_2D_Traj}
  \end{minipage}
  \hfill
  \begin{minipage}{0.49\textwidth}
    \includegraphics[width=\textwidth]{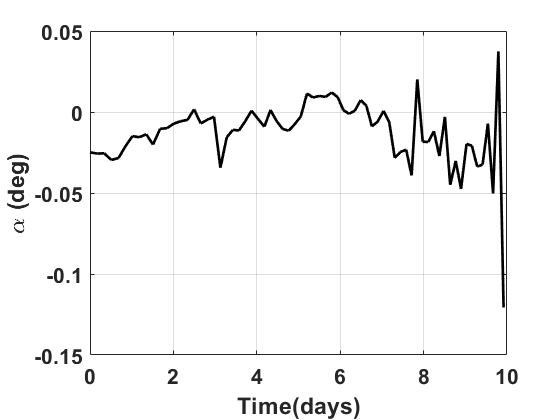}
    \caption{Case 2: History of $\alpha$}
    \label{fig:Case2_LEO_MEO_2D_Alpha}
  \end{minipage}
\end{figure}

The transfer maneuver from a low Earth circular orbit to a medium Earth circular orbit is considered here. The initial and final altitudes are $h_i = 300 \, km$ and $h_f = 2000 \, km$ respectively. The parameters of the spacecraft are defined as follows:
specific impulse $Isp$ is $3300 \, s$, initial spacecraft mass is $95 \,kg$, thrust-to-weight ratio = $10^{-4}$ and maximum
thrust limit is $0.0920 \, N$. 
\begin{table}[hbt!]
\caption{\label{Case2_DP} Case 2: Design Parameters}
\centering
\begin{tabular}{ccc}
\hline
\hline
        Parameter & Initial Orbit & Target Orbit \\ 
        \hline
        Semi-major axis, $a$ (km) & 6678.18 & 8378.18 \\ 
        Eccentricity, $e$ (km) & 0 & 0 \\ 
        Inclination, $i$ (deg) & 0 & 0 \\
        Argument of Periapsis, $\omega$ (deg) & 0 & 0 \\
        RAAN, $\Omega$ (deg) & 0 & 0 \\
        True Anomaly, $\theta_f$ & 0$^{\circ}$ & 90$^{\circ}$ \\
\hline
\hline
\end{tabular}
\end{table}
The boundary conditions(BC) of the initial and target orbit at initial time is tabulated in Table~\ref{Case2_DP}. The corresponding trajectory solution is shown in Fig~\ref{fig:Case2_LEO_MEO_2D_Traj}.
The entire trajectory was divided into 68 segments each corresponding to an altitude increase of 25 km. The spacecraft performs 268 revolutions around the Earth in 9.91 days to achieve rendezvous condition with a secondary spacecraft in the MEO. To achieve this maneuver, the spacecraft consumes 2.4 kg of propellant with an associated $\Delta V$ = 0.8406 km/s for this transfer. Figure~\ref{fig:Case2_LEO_MEO_2D_Alpha} shows the history of the thrust steering angle $\alpha$. The performance characteristics of the obtained solution are tabulated in Table~\ref{tab:Case2_Performance_Characteristics}
\begin{figure}[h!]
  \centering
  \begin{minipage}{0.49\textwidth}
    \includegraphics[width=\textwidth]{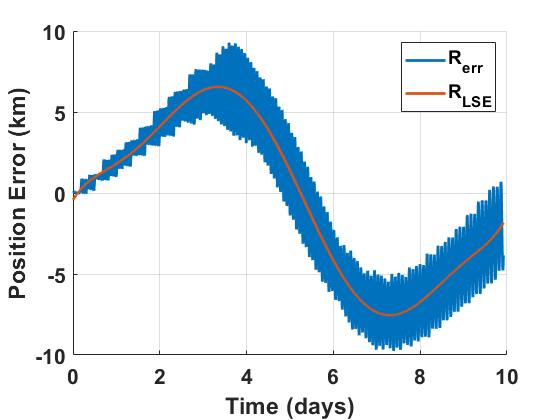}
    \caption{Case 2: Positional Error}
    \label{fig:Case2_LEO_MEO_2D_r_err}
  \end{minipage}
  \hfill
  \begin{minipage}{0.49\textwidth}
    \includegraphics[width=\textwidth]{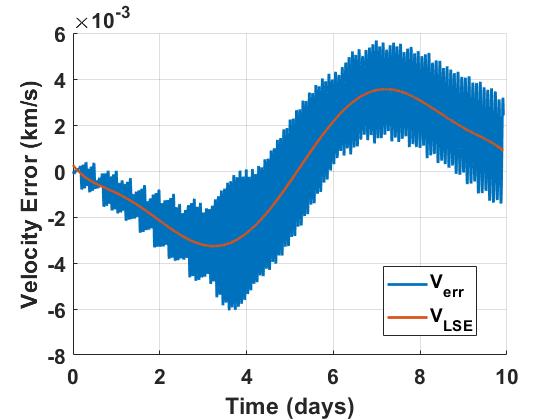}
    \caption{Case 2: Velocity Error}
    \label{fig:Case2_LEO_MEO_2D_v_err}
  \end{minipage}
\end{figure}

\begin{table}[h!]
	\fontsize{10}{10}\selectfont
    \caption{Case 2: Performance Characteristics of the Semi-Analytic Solution}
   \label{tab:Case2_Performance_Characteristics}
        \centering 
   \begin{tabular}{c | c } 
      \hline 
        Parameters & Value \\ 
            \hline
        No. of segments & 68\\
        No. of revolutions & 268 \\
        Time of Flight (days) & 9.91 \\ 
        $\Delta V$ (km/s) & 0.8406 \\
        Propellant Consumption (kg) & 2.4 \\
      \hline
   \end{tabular}
\end{table}

Figures~\ref{fig:Case2_LEO_MEO_2D_r_err} and ~\ref{fig:Case2_LEO_MEO_2D_v_err} illustrate the positional and velocity error between the semi-analytic solution and two-body numerical propagation of the spacecraft. Solid black curves are used to highlight the least square fit of the positional and velocity error. The final states of the spacecraft achieved by the semi-analytic method and numerical propagation are compared in Table~\ref{tab:Case2_CW_NP_Comp} to demonstrate the level of error between them. 
\begin{table}[hbt!]
\caption{\label{tab:Case2_CW_NP_Comp} Case 2: Semi-Analytic Solution vs Numerical Propagation}
\centering
\begin{tabular}{ccc}
\hline
\hline
        Target orbit Parameters & Semi-Analytic Solution & Numerical Propagation \\ 
        \hline
        Semi-major axis, $a$ (km) & 8378.18 & 8375.95 \\ 
        Eccentricity, $e$ (km) & 0.00002 & 0.0003 \\ 
        Inclination, $i$ (deg) & 0 & 0 \\
        Argument of Periapsis, $\omega$ (deg) & 0 & 0 \\
        RAAN, $\Omega$ (deg) & 0 & 0 \\
        True Anomaly, $\theta_f$ & 60$^{\circ}$ & 50$^{\circ}$ \\
\hline
\hline
\end{tabular}
\end{table}


\subsection{Case 3: 3D Transfer LEO-to-GEO (Orbit Insertion)}\label{Case2} 

A three dimensional transfer from an initial circular low Earth orbit (LEO) at $h_i = 300 \,$ km altitude and $i_i = 28.5^{\circ}$  inclination to a final circular Geostationary Orbit (GEO) of altitude $h_f = 35,785\,km$ and $i_f = 0^{\circ}$ inclination is considered in this section. This LEO-to-GEO orbit insertion problem is of practical importance considering the remote sensing, weather forecasting, Earth observation applications and hence extensively studied\cite{taheri2012shape}\cite{kluever1998simple}\cite{shepard2020preliminary} in literature. The initial orbit was so chosen since the Kennedy Space Center is at $28.5^{\circ}$ North latitude and hence cheaper to launch spacecrafts into a parking orbit on the same latitude. In this study, the parameters of the spacecraft are defined as follows: specific impulse, $Isp = 2800 \,s$, initial spacecraft mass, $m_0 = 95 \,kg$, thrust-to-weight ratio = $10^{-4}$ and maximum thrust limit is $0.0920 \,N$.
\begin{figure}[h!]
  \centering
  \begin{minipage}{0.49\textwidth}
    \includegraphics[width=\textwidth]{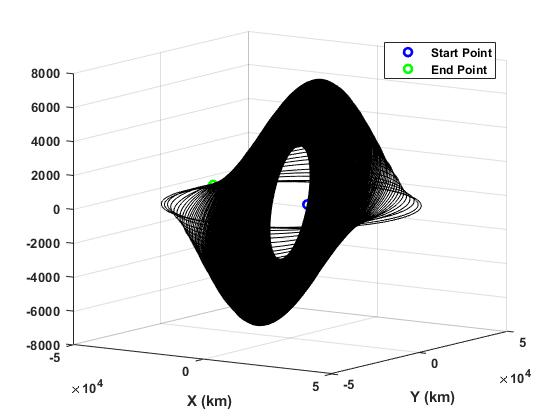}
    \caption{Case 3: Trajectory (ECI Frame)}
    \label{fig:Case3_LEO_GEO_3D_Traj}
  \end{minipage}
  \hfill
  \begin{minipage}{0.49\textwidth}
    \includegraphics[width=\textwidth]{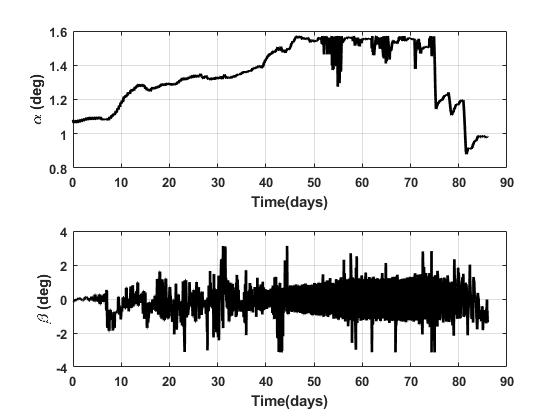}
    \caption{Case 3: History of $\alpha$ and $\beta$}
    \label{fig:Case3_LEO_GEO_3D_Alpha}
  \end{minipage}
\end{figure}
The entire trajectory design problem was divided into 3550 segments, roughly corresponding to an altitude change of $10 \, km$ per segment. The trajectory solution and the history of the thrust steering angles are shown in Figure~\ref{fig:Case3_LEO_GEO_3D_Traj} and ~\ref{fig:Case3_LEO_GEO_3D_Alpha}. The out-of-plane thrust component is used to reduce the initial inclination set by the parking orbit, while the in-plane thrust component raises the orbit semi-major axis to the desired geostationary level. The obtained solution performs $N_{rev} = 512$ revolutions and has a total flight time of $86.01 \,days$. The total $\Delta V$ for the transfer is $7.29 \, km/s$ and the spacecraft consumes $21.86 \, kg$ of propellant. 

The least square fit of the errors are represented using solid black curves. The performance characteristics of the obtained solution are tabulated in Table~\ref{tab:Case3_Performance_Characteristics}.
\begin{table}[h!]
	\fontsize{10}{10}\selectfont
    \caption{Case 3: Performance Characteristics of the Semi-Analytic Solution}
   \label{tab:Case3_Performance_Characteristics}
        \centering 
   \begin{tabular}{c | c } 
      \hline 
        Parameters & Value \\ 
            \hline
        No. of segments & 3550\\
        No. of revolutions & 512 \\
        Time of Flight (days) & 86.01 \\ 
        $\Delta V$ (km/s) & 7.29 \\
        Propellant Consumption (kg) & 21.86\\
      \hline
   \end{tabular}
\end{table}

\subsection{Case 4: 3D Transfer Earth-to-Mars (Rendezvous)}\label{Case2} 

\begin{figure}[h!]
  \centering
  \begin{minipage}{0.49\textwidth}
    \includegraphics[width=\textwidth]{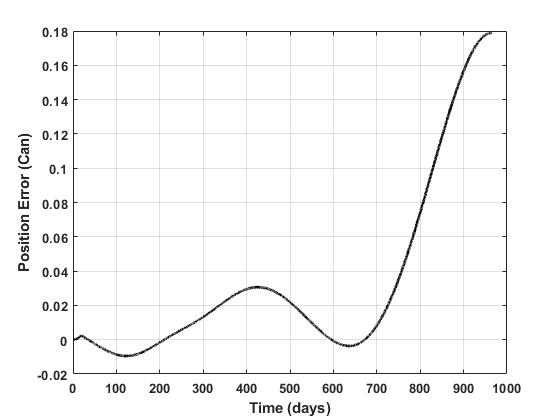}
    \caption{Case 4: Positional Error}
    \label{fig:Case4_Earth_Mars_3D_r_err}
  \end{minipage}
  \hfill
  \begin{minipage}{0.49\textwidth}
    \includegraphics[width=\textwidth]{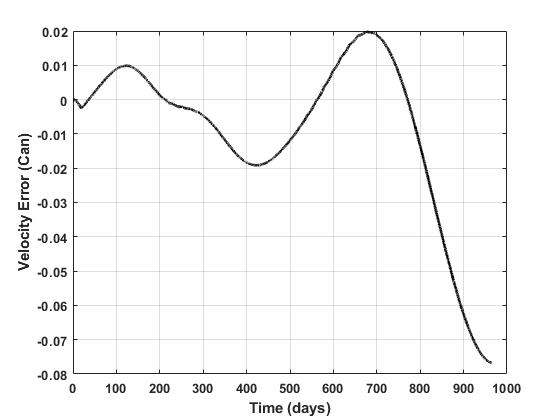}
    \caption{Case 4: Velocity Error}
    \label{fig:Case2_Earth_Mars_3D_v_err}
  \end{minipage}
\end{figure}

A low thrust Earth-Mars transfer problem is studied in this section. The spacecraft is expected to rendezvous with Mars starting from Earth given the low thrust limitation. This problem is of prime interest to the aerospace community as the global community of scientists and engineers hope to colonize Mars within the next few decades. Low thrust trajectories would be key to facilitate regular transportation of supplies and cargo between Earth and Mars. In this case study, the date of launch of the spacecraft departing from Earth is assumed to be fixed and arbitrarily chosen to be the 20th of July 2023 marking 54 years since the first Moon landing. The orbital elements of Earth and Mars on this date is shown in Table~\ref{Case4_DP}.
\begin{table}[hbt!]
\caption{\label{Case4_DP} Case 4: Design Parameters}
\centering
\begin{tabular}{ccc}
\hline
\hline
        Parameter & Earth Orbit & Mars Orbit \\ 
        \hline
        Semi-major axis (AU) & 1 & 1.52366231 \\ 
        Eccentricity & 0.01671022 & 0.09341233 \\ 
        Inclination (deg) & 0.00005 & 1.85061 \\
        Longitude of Perihelion (deg) & 102.94719 & 336.04084 \\
        Longitude of ascending node (deg) & -11.26064 & 49.57854 \\
        True Anomaly (deg) & 194.72$^{\circ}$ & 201.99$^{\circ}$ \\
\hline
\hline
\end{tabular}
\end{table}
The spacecraft propulsion system characteristics for this study are detailed as follows: specific impulse, $Isp = 2800 \, s$, initial
spacecraft mass, $m_0 = 1000 \, kg$, thrust-to-weight ratio = $10^{-6}$ and maximum thrust limit = $0.098 \,N$.
\begin{figure}[h!]
  \centering
  \begin{minipage}{0.49\textwidth}
    \includegraphics[width=\textwidth]{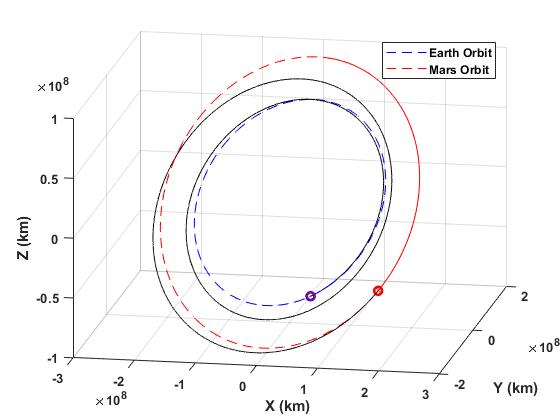}
    \caption{Case 4: Trajectory (ICRF)}
    \label{fig:Case4_Earth_Mars_3D_Traj}
  \end{minipage}
  \hfill
  \begin{minipage}{0.49\textwidth}
    \includegraphics[width=\textwidth]{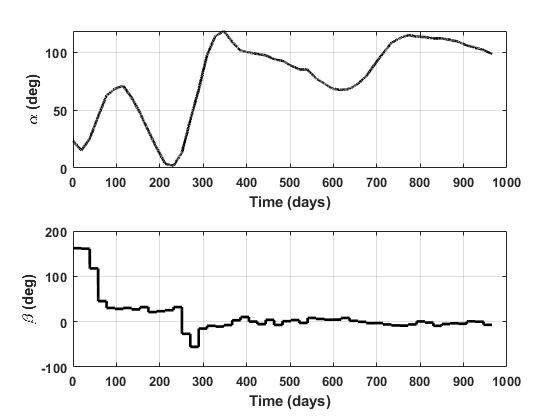}
    \caption{Case 4: History of $\alpha$ and $\beta$}
    \label{fig:Case4_Earth_Mars_3D_Alpha}
  \end{minipage}
\end{figure}
The trajectory design process is carried out by dividing the entire trajectory into 50 segments. All the calculations for this case study are carried out using canonical units such that one distance
unit (DU) is equal to the radius of the Earth's orbit around the Sun ($1.496e+8$ km). The trajectory solution obtained by the semi-analytic method is shown in Fig.~\ref{fig:Case4_Earth_Mars_3D_Traj}. It is important to note that the algorithm simply returns a feasible rendezvous trajectory and no optimal solution is claimed by this method. The history of the thrust steering angles are shown in Fig.~\ref{fig:Case4_Earth_Mars_3D_Alpha}. Note that the states (position and velocity) of Earth and Mars for this problem was obtained from the DE430\cite{folkner2014planetary} Ephemeris data found on the JPL website. The spacecraft takes $965.33 \,$days to rendezvous with Mars. In other words, the spacecraft arrives at Mars on 11th March 2026. To perform this transfer, the spacecraft consumes $261.87 \, kg$ of propellant, performing 2 revolutions around the Earth. The total $\Delta V$ of this transfer is $8.34 \, km/s$. The performance characteristics of the trajectory solution are tabulated in Table~\ref{tab:Case4_Performance_Characteristics}.
\begin{table}[h!]
	\fontsize{10}{10}\selectfont
    \caption{Case 4: Performance Characteristics of the Semi-Analytic Solution}
   \label{tab:Case4_Performance_Characteristics}
        \centering 
   \begin{tabular}{c | c } 
      \hline 
        Parameters & Value \\ 
            \hline
        No. of segments & 50\\
        No. of revolutions & 2 \\
        Time of Flight (days) & 965.33 \\ 
        $\Delta V$ (km/s) & 8.34 \\
        Propellant Consumption (kg) & 261.87\\
      \hline
   \end{tabular}
\end{table}
Figures~\ref{fig:Case4_Earth_Mars_3D_r_err} and \ref{fig:Case2_Earth_Mars_3D_v_err} highlight the positional and velocity error between the semi-analytic solution and two-body numerical propagation of the spacecraft.
\section{Conclusion}

The work developed in this paper presents a semi-analytic approach for the generation of initial guess guidance trajectories for low-thrust spacecrafts. A modification on the Hill's equations is provided as a means to approximate the states of the spacecraft acted upon by a constant low-thrust acceleration. Numerical results demonstrates the flexibility of the algorithm in generating three-dimensional rendezvous, orbit insertion and orbit raising trajectory solutions.
\FloatBarrier
\section*{Acknowledgement}
This paper is based upon work supported by NASA, Grant Number 80NSSC19K1642

\begin{landscape}
\section{APPENDIX-A}
Analytic expression for the in-plane motion of the spacecraft.
\begin{dmath}
    x(t) = \frac{2\,k^3 \,{\dot{y} }_0 -4\,k\,n^3 \,x_0 +4\,k^3 \,n\,x_0 -2\,k\,n^2 \,{\dot{y} }_0 -2\,k^3 \,{\dot{y} }_0 \,\mathrm{cos}\left(n\,t\right)}{k\,n\,{\left(k^2 -n^2 \right)}} +\frac{k^3 \,{\dot{x} }_0 \,\mathrm{sin}\left(n\,t\right)-2\,a\,k^2 \,\mathrm{cos}\left(\alpha_0 \right)\,\mathrm{cos}\left(\beta \right)+2\,a\,n^2 \,\mathrm{cos}\left(\alpha_0 \right)\,\mathrm{cos}\left(\beta \right)}{k\,n\,{\left(k^2 -n^2 \right)}} +\frac{3\,k\,n^3 \,x_0 \,\mathrm{cos}\left(n\,t\right)-3\,k^3 \,n\,x_0 \,\mathrm{cos}\left(n\,t\right)+2\,k\,n^2 \,{\dot{y} }_0 \,\mathrm{cos}\left(n\,t\right)}{k\,n\,{\left(k^2 -n^2 \right)}} +\frac{-k\,n^2 \,{\dot{x} }_0 \,\mathrm{sin}\left(n\,t\right) + 2\,a\,k^2 \,\mathrm{cos}\left(n\,t\right)\,\mathrm{cos}\left(\alpha_0 \right)\,\mathrm{cos}\left(\beta \right)-2\,a\,n^2 \,\mathrm{cos}\left(k\,t\right)\,\mathrm{cos}\left(\alpha_0 \right)\,\mathrm{cos}\left(\beta \right)}{k\,n\,{\left(k^2 -n^2 \right)}} +\frac{a\,k^2 \,\mathrm{sin}\left(n\,t\right)\,\mathrm{cos}\left(\beta \right)\,\mathrm{sin}\left(\alpha_0 \right)-2\,a\,n^2 \,\mathrm{sin}\left(k\,t\right)\,\mathrm{cos}\left(\beta \right)\,\mathrm{sin}\left(\alpha_0 \right)-a\,k\,n\,\mathrm{cos}\left(k\,t\right)\,\mathrm{cos}\left(\alpha_0 \right)\,\mathrm{cos}\left(\beta \right)}{k\,n\,{\left(k^2 -n^2 \right)}} + \frac{+a\,k\,n\,\mathrm{cos}\left(n\,t\right)\,\mathrm{cos}\left(\alpha_0 \right)\,\mathrm{cos}\left(\beta \right) -a\,k\,n\,\mathrm{sin}\left(k\,t\right)\,\mathrm{cos}\left(\beta \right)\,\mathrm{sin}\left(\alpha_0 \right)+2\,a\,k\,n\,\mathrm{sin}\left(n\,t\right)\,\mathrm{cos}\left(\beta \right)\,\mathrm{sin}\left(\alpha_0 \right)}{k\,n\,{\left(k^2 -n^2 \right)}}
\end{dmath}
\begin{dmath}
\dot{x}(t) = \frac{k^2 \,{\dot{x} }_0 \,\mathrm{cos}\left(n\,t\right)-n^2 \,{\dot{x} }_0 \,\mathrm{cos}\left(n\,t\right)+2\,k^2 \,{\dot{y} }_0 \,\mathrm{sin}\left(n\,t\right)}{k^2 -n^2 } + \frac{-3\,n^3 \,x_0 \,\mathrm{sin}\left(n\,t\right)-2\,n^2 \,{\dot{y} }_0 \,\mathrm{sin}\left(n\,t\right)+3\,k^2 \,n\,x_0 \,\mathrm{sin}\left(n\,t\right)}{k^2 -n^2 } + \frac{-a\,k\,\mathrm{cos}\left(k\,t\right)\,\mathrm{cos}\left(\beta \right)\,\mathrm{sin}\left(\alpha_0 \right)+a\,k\,\mathrm{sin}\left(k\,t\right)\,\mathrm{cos}\left(\alpha_0 \right)\,\mathrm{cos}\left(\beta \right)+a\,k\,\mathrm{cos}\left(n\,t\right)\,\mathrm{cos}\left(\beta \right)\,\mathrm{sin}\left(\alpha_0 \right)}{k^2 -n^2 } + \frac{-2\,a\,k\,\mathrm{sin}\left(n\,t\right)\,\mathrm{cos}\left(\alpha_0 \right)\,\mathrm{cos}\left(\beta \right)-2\,a\,n\,\mathrm{cos}\left(k\,t\right)\,\mathrm{cos}\left(\beta \right)\,\mathrm{sin}\left(\alpha_0 \right)+2\,a\,n\,\mathrm{sin}\left(k\,t\right)\,\mathrm{cos}\left(\alpha_0 \right)\,\mathrm{cos}\left(\beta \right)}{k^2 -n^2 } +\frac{2\,a\,n\,\mathrm{cos}\left(n\,t\right)\,\mathrm{cos}\left(\beta \right)\,\mathrm{sin}\left(\alpha_0 \right)-a\,n\,\mathrm{sin}\left(n\,t\right)\,\mathrm{cos}\left(\alpha_0 \right)\,\mathrm{cos}\left(\beta \right)}{k^2 -n^2 }
\end{dmath}
\end{landscape}
\begin{landscape}
\begin{dmath}
y(t) = y_0 + \frac{ t\,{\left(3\,k^2 -3\,n^2 \right)}\,{\left(k\,{\dot{y} }_0 -a\,\mathrm{cos}\left(\alpha_0 \right)\,\mathrm{cos}\left(\beta \right)+2\,k\,n\,x_0 \right)}}{k\,\left(n^2 -k^2\right) } +\frac{2\,k\,\mathrm{sin}\left(n\,t\right)\,{\left(-3\,x_0 \,k^2 \,n-2\,{\dot{y} }_0 \,k^2 +2\,a\,\mathrm{cos}\left(\alpha_0 \right)\,\mathrm{cos}\left(\beta \right)\,k+3\,x_0 \,n^3 +2\,{\dot{y} }_0 \,n^2 +a\,\mathrm{cos}\left(\alpha_0 \right)\,\mathrm{cos}\left(\beta \right)\,n\right)}}{k\,\left(n^3 -k^2\right)} -\frac{2\,k\,\mathrm{cos}\left(n\,t\right)\,{\left({\dot{x} }_0 \,k^2 +a\,\mathrm{cos}\left(\beta \right)\,\mathrm{sin}\left(\alpha_0 \right)\,k-{\dot{x} }_0 \,n^2 +2\,a\,\mathrm{cos}\left(\beta \right)\,\mathrm{sin}\left(\alpha_0 \right)\,n\right)}}{k\,\left(n^3 -k^2\right)} +\frac{a\,\mathrm{sin}\left(\alpha_0 -k\,t\right)\,\mathrm{cos}\left(\beta \right)\,{\left(k^2 -n^2 \right)}}{k^2\left(n^2 -k \right)} +\frac{2\,a\,n\,\mathrm{cos}\left(k\,t\right)\,\mathrm{cos}\left(\beta \right)\,\mathrm{sin}\left(\alpha_0 \right)\,{\left(k+2\,n\right)}}{k^2\,\left(n^2 -k \right)} -\frac{2\,a\,n\,\mathrm{sin}\left(k\,t\right)\,\mathrm{cos}\left(\alpha_0 \right)\,\mathrm{cos}\left(\beta \right)\,{\left(k+2\,n\right)}}{k^2\,\left(n^2 -k \right)} -\frac{2\,{\dot{x} }_0 \,k^2 +2\,a\,\mathrm{cos}\left(\beta \right)\,\mathrm{sin}\left(\alpha_0 \right)\,k+3\,a\,n\,\mathrm{cos}\left(\beta \right)\,\mathrm{sin}\left(\alpha_0 \right)}{k^2 \,n}
\end{dmath}

\begin{dmath}
\dot{y}(t) = {\dot{y} }_0 +\frac{a\,\mathrm{cos}\left(\alpha_0 -k\,t\right)\,\mathrm{cos}\left(\beta \right)\,\sigma_1 }{k\,\left(k^2 - n^2 \right)} -\frac{\mathrm{cos}\left(n\,t\right)\,{\left(-6\,x_0 \,k^2 \,n-4\,{\dot{y} }_0 \,k^2 +4\,a\,\mathrm{cos}\left(\alpha_0 \right)\,\mathrm{cos}\left(\beta \right)\,k+6\,x_0 \,n^3 +4\,{\dot{y} }_0 \,n^2 +2\,a\,\mathrm{cos}\left(\alpha_0 \right)\,\mathrm{cos}\left(\beta \right)\,n\right)}}{k^2 - n^2} -\frac{\mathrm{sin}\left(n\,t\right)\,{\left(2\,{\dot{x} }_0 \,k^2 +2\,a\,\mathrm{cos}\left(\beta \right)\,\mathrm{sin}\left(\alpha_0 \right)\,k-2\,{\dot{x} }_0 \,n^2 +4\,a\,\mathrm{cos}\left(\beta \right)\,\mathrm{sin}\left(\alpha_0 \right)\,n\right)}}{k^2 - n^2} +\frac{2\,a\,n\,\mathrm{cos}\left(k\,t\right)\,\mathrm{cos}\left(\alpha_0 \right)\,\mathrm{cos}\left(\beta \right)\,{\left(k+2\,n\right)}}{k\,\left(k^2 - n^2 \right)} +\frac{2\,a\,n\,\mathrm{sin}\left(k\,t\right)\,\mathrm{cos}\left(\beta \right)\,\mathrm{sin}\left(\alpha_0 \right)\,{\left(k+2\,n\right)}}{k\,\left(k^2 - n^2 \right)} -\frac{4\,k\,{\dot{y} }_0 -3\,a\,\mathrm{cos}\left(\alpha_0 \right)\,\mathrm{cos}\left(\beta \right)+6\,k\,n\,x_0 }{k}
\end{dmath}
\end{landscape}

\end{document}